\documentclass{elsart}
\usepackage{amssymb}

\begin{document}

\begin{frontmatter}

\title{Statistical duality of the Laplace distribution}

\author[1]{E.A. Barkova}
\author[2]{S.I. Bityukov\corauthref{cor1}}
\corauth[cor1]{Serguei.Bitioukov@cern.ch}
\author[1]{V.A. Taperechkina}
\address[1]{Moscow State Academy of Instrument Engineering and 
            Computer Science, Moscow, Russia}
\address[2]{Institute for high energy physics, 142281 Protvino, Russia}

\begin{abstract}
The statistical duality of distributions is a powerful tool
for statistical inferences. In the paper the statistical duality
of Laplace distribution is discussed. As shown the confidence density 
of the parameter of this distribution is uniquely determined.
\end{abstract}
  

\end{frontmatter}

\section{Introduction}
\label{Int}

In ref.\cite{Bit2004} is introduced the notion of 
statistical duality and is shown that several pairs of distributions 
(Poisson and Gamma, normal and normal, Cauchy and Cauchy) are 
statistically dual distributions. These distributions  allow to exchange the 
parameter and the random variable, conserving the same formula for the 
distribution of probabilities. The interrelation between the statistically 
dual distributions and conjugate families is considered
in ref.\cite{Bit2005}. It allows to use the statistical duality for 
estimation of the distribution parameter.

In the paper we show that Laplace distributions are statistically 
dual distributions. 

\section{Statistical duality and estimation of the parameter of distribution} 
\label{EPD}

The probability density of the Laplace distribution is

\begin{equation}
\displaystyle L(x;a,b) = 
\displaystyle \frac{1}{2b} \displaystyle  e^{-\frac{|x-a|}{b}},
\end{equation}

\noindent
where $x$ is a real variable, $a$ and $b>0$ are real parameters.
Here we can exchange the parameter $a$ and variable $x$ with conserving 
of the same formula for new density. This new density can be named 
as a confidence density of the parameter $a$~(for example,~\cite{Efr1998}) 

\begin{equation}
\displaystyle \tilde L(a;x,b) = 
\displaystyle \frac{1}{2b} \displaystyle  e^{-\frac{|x-a|}{b}},
\end{equation}

\noindent
where $a$ is a real variable, $x$ and $b>0$ are real parameters.

\noindent
By this means the Laplace distributions are a statistically self-dual 
distributions~\cite{Bit2004}. 

Let $\hat x$ be a result of single observation of random variable $x$.
If we define the confidence interval 
$(a_1,a_2)$~(see, also, \cite{Bit2004,Bit2000}) 
as~\footnote{Note, in this case the given definition coincides
with definition of fiducial interval~\cite{Fis1930}.} 

\begin{equation}
P(a_1 \le a \le a_2|\hat x) = 
P(x \le \hat x | a_1) - P(x \le \hat x | a_2),
\end{equation}

\noindent
where $a_1 \le a_2$ are real values and $P(x \le \hat x | a) = 
\displaystyle \int_{-\infty}^{\hat x}{L(x;a,b)dx}$,
we can reconstruct the confidence density of the 
parameter $a$. This is due to the fact that the identity

\begin{equation}
\displaystyle \int_{\hat x}^{\infty}{L(x;a_1,b)dx} + 
\int_{a_1}^{a_2}{\tilde L(a;\hat x,b)da} +
\displaystyle \int_{-\infty}^{\hat x}{L(x;a_2,b)dx} = 1
\end{equation}

\noindent
takes place for any real values $b>0$ and $a_1 \le a_2$. 

Let us suppose that $\tilde L(a;\hat x,b)$ is the confidence density 
of parameter of the Laplace distribution
if observed value of random variable $x$ is equal to $\hat x$.
It is a conditional probability density. As it follows from formulae 
(Eqs.2,4), the $\tilde L(a;\hat x,b)$ is the density of Laplace distribution 
by definition. 

On the other hand: 
if $\tilde L(a;\hat x,b)$ is not equal to this confidence density 
and the confidence density of the Laplace parameter is the other 
function $h(a;\hat x,b)$ then there takes place another identity

\begin{equation}
\displaystyle \int_{\hat x}^{\infty}{L(x;a_1,b)dx} + 
\int_{a_1}^{a_2}{h(a;\hat x,b)da} +
\displaystyle \int_{-\infty}^{\hat x}{L(x;a_2,b)dx} = 1
\end{equation}

This identity is correct for any real $a_1 \le a_2$  and 
$\hat x$ too.
The first and third terms in the left part of this identity determine the 
boundary conditions on the confidence interval.

If we subtract Eq.5 from Eq.4 then we have

\begin{equation}
\int_{a_1}^{a_2}{\tilde L(a;\hat x,b) - h(a;\hat x,b)da}  = 0.
\end{equation}

We can choose the $a_1$ and $a_2$ by the arbitrary way. Let us make this 
choice so that $\tilde L(a;\hat x,b)$ is not equal $h(a;\hat x,b)$ 
in the interval
$(a_1,a_2)$ and, for example, $\tilde L(a;\hat x,b) > h(a;\hat x,b)$ and
$a_2 > a_1$. In this case we have 

\begin{equation}
\int_{a_1}^{a_2}{\tilde L(a;\hat x,b) - h(a;\hat x,b)da}  > 0.
\end{equation}

\noindent
and we have contradiction. Hence $\tilde L(a;\hat x,b)=h(a;\hat x,b)$ 
everywhere except, may be, a finite set of points.  

As a consequence, the reconstruction
of the confidence density of the Laplace distribution parameter $a$ is an 
unique, i.e. the confidence density of the parameter $a$ is the probability 
density $\tilde L(a;\hat x,b)$.

So, the statistical duality allows to connect the estimation
of the parameter with the measurement of the random variable of the
Laplace distribution.

\section{Conclusion}

In the paper we show that the Laplace distribution is
a statistically self-dual distributions. The confidence density of the
Laplace distribution parameter is determined in case of single observation
$\hat x$ of random variable $x$.
It allows to construct the confidence intervals for parameter $a$ of 
the distribution by the easy way.

\subsubsection*{Acknowledgments}

The authors are grateful to V.A.~Kachanov, N.V.~Krasnikov 
and V.F.~Obraztsov for the interest and useful comments, 
S.S.~Bityukov, Yu.P.~Gouz, C.~Wulz 
for fruitful discussions and E.A.~Medvedeva for help in preparing
the paper. This work has been particularly supported 
by grants RFBR 04-01-97227 and RFBR 04-02-16020.

\end{document}